\documentclass[11pt]{article}

\usepackage{amsthm,amsmath,amssymb,ifthen,macros,pstricks,pst-node}
\usepackage[breaklinks]{hyperref}
\usepackage[round,sort&compress]{natbib}
\usepackage[top=3cm,bottom=3cm,left=3.75cm,right=3.75cm]{geometry}

\usepackage[sort&compress,round]{natbib}
\usepackage[breaklinks]{hyperref}

\DeclareMathOperator{\Ew}{E}       
\DeclareMathOperator{\Var}{Var}       
\DeclareMathOperator{\vecc}{vec}       
\DeclareMathOperator{\lsum}{\textstyle\sum}
\numberwithin{equation}{section}

\begin{document}

\title{Computing Maximum Likelihood Estimates in Recursive Linear Models
  with Correlated Errors}

\author{Mathias Drton
  \\
  University of Chicago\\
  {\em e-mail:} drton@uchicago.edu\\
  \and
  Michael Eichler\\
  University of Maastricht\\
  {\em e-mail:} m.eichler@ke.unimaas.nl\\
  \and
  Thomas S.~Richardson
  \\
  University of Washington \\
  {\em e-mail:} tsr@stat.washington.edu
} 

\date{} \maketitle 

\begin{abstract}
  In recursive linear models, the multivariate normal joint distribution of all variables exhibits a dependence structure induced by a recursive (or acyclic) system of linear structural equations.  These linear models have a long tradition and appear in seemingly unrelated regressions, structural equation modelling, and approaches to causal inference.  They are also related to Gaussian graphical models via a classical  representation known as a path diagram.  Despite the models' long   history, a number of problems remain open.  In this paper, we address the   problem of computing maximum likelihood estimates in the subclass of   `bow-free' recursive linear models. The term `bow-free' refers to the   condition that the errors for variables $i$ and $j$ be uncorrelated if   variable $i$ occurs in the structural equation for variable $j$.  We   introduce a new algorithm, termed Residual Iterative Conditional Fitting  (RICF), that can be implemented using only least squares computations.  In contrast to existing algorithms, RICF has clear convergence properties and finds parameter estimates in closed form whenever possible.\\

\end{abstract} 


\renewcommand{\theenumi}{\roman{enumi}}
\renewcommand{\labelenumi}{(\theenumi)}
\newcommand{\mss}{\scriptstyle}
\newcommand{\I}[2]{ {#1}{\ind}{#2}} 
\newcommand{\CI}[3]{ {#1}{\ind}{#2}\,|\, {#3} }
\newcommand{\db}{\text{\rm db}}       
\newcommand{\idmat}{I}  
\newcommand{\veps}{\varepsilon} 

\section{Introduction}
\label{sec:introduction}

A system of linear structural equations determines a linear model for a set
of variables by dictating that, up to a random error term, each variable is
equal to a linear combination of some of the remaining variables.
Traditionally the errors are assumed to have a centered joint multivariate
normal distribution.  Presenting a formalism for simultaneously
representing causal and statistical hypotheses
\citep{pearl:causalbk,spirtes:2000}, these normal linear models, which are
often called {\em structural equation models\/}, are widely used in the
social sciences \citep{bollen:89} and many other contexts.

In seminal work, \citet{wright:correlation:1921, wright:method:1934}
introduced {\em path diagrams}, which are very useful graphical
representations of structural equations.  A path diagram is a graph with
one vertex for each variable and directed and/or bi-directed edges.  A
directed edge $i\rightarrow j$ indicates that variable $i$ appears as
covariate in the equation for variable $j$.  The directed edges are thus in
correspondence with the {\em path coefficients}, that is, the coefficients
appearing in the linear structural equations.  A bi-directed edge
$i\leftrightarrow j$ indicates correlation between the errors in the
equations for variables $i$ and $j$. Graphs of this kind are also considered
by \citet{conf:aaai:ShpitserP06}, who  refer to them as
 recursive semi-Markovian causal models.

\subsection{A Motivating Example}
\label{sec:motivating-example}

We will motivate the considered normal linear models with the following
example, which is adapted from a more complex longitudinal study considered
in \cite{obesity:robins:2008}.

Consider a two-phase sequential intervention study examining the effect of
exercise and diet on diabetes. In the first phase patients are randomly
assigned to a number of hours of exercise per week (Ex) drawn from a
log-normal distribution. At the end of this phase blood pressure (BP)
levels are measured.  In the second phase patients are randomly assigned to
a strict calorie controlled diet that produces a change in body-mass index
($\Delta$BMI).  The assigned change in BMI, though still randomized, is
drawn, by design, from a normal distribution with mean depending linearly
on $\log(\hbox{Ex})$ and BP. The dependence here is due to practical and
ethical considerations. Finally at the end of the second phase,
triglyceride levels ($Y$) indicating diabetic status are measured.

A question of interest is whether or not there is an effect of Ex on the
outcome $Y$ that is not mediated through the dependence of $\Delta$BMI on
Ex and BP.  In other words, if there had been no ethical or practical
restrictions, and the assignment ($\Delta$BMI) in the second phase was
completely randomized and thus independent of BP and Ex, would there still
be any dependence between Ex and $Y$? Note that due to underlying
confounding factors such as life history and genetic background, we would
expect to observe dependence between BP and $Y$ even if the null hypothesis
of no effect of Ex on $Y$ was true and the second treatment ($\Delta$BMI)
was completely randomized. 

\begin{figure}[t]
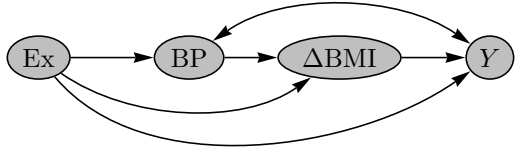

\begin{center}
    \small
    \vspace{2cm}\hspace{-5cm}
     \psset{linewidth=0.6pt,arrowscale=1.5 2,arrowinset=0.1}          
    \psset{fillcolor=lightgray, fillstyle=solid}          
    \newcommand{\myNode}[2]{\ovalnode{#1}{#2}}
    \rput(-1, 1){\myNode{1}{Ex}}
    \rput(1, 1){\myNode{2}{BP}} 
    \rput(3, 1){\myNode{3}{$\Delta$BMI}} 
    \rput(5, 1){\myNode{4}{$Y$}} 
    \psset{fillcolor=lightgray, fillstyle=none}          
    \ncline{->}{1}{2}
    \ncline{->}{2}{3}
    \ncline{->}{3}{4}
 \nccurve[angleA=-35,angleB=215]{->}{1}{3}
  \nccurve[angleA=-50,angleB=215]{->}{1}{4}
    \nccurve[angleA=35,angleB=145]{<->}{2}{4}
   \end{center}
\caption{\label{fig:0robins} Path diagram illustrating a two-phase 
  trial with two treatments (Ex and $\Delta$BMI) and two responses (BP and
  $Y$). Ex is randomly assigned, $\Delta$BMI is randomized conditional on
  BP and Ex. The bi-directed edge indicates possible dependence due to
  unmeasured factors (genetic or environmental).}
\end{figure}

Our model consists of two pieces.  First, the design of the study dictates
that
\begin{align}
  \log (\hbox{Ex}) &= \alpha_0 + \varepsilon_{\hbox{\tiny Ex}},\label{eq:ex}\\
  \Delta \hbox{BMI} &= \gamma_0 + \gamma_1 \log(\hbox{Ex}) +\gamma_2\,
  \hbox{BP} +   \varepsilon_{\Delta \hbox{\tiny BMI}},\label{eq:bmi}
\end{align}
where $\varepsilon_{\hbox{\tiny Ex}}\sim \mathcal{N}(0,\sigma^2_{\hbox{\tiny Ex}})$
and $\varepsilon_{\Delta \hbox{\tiny BMI}}\sim
\mathcal{N}(0,\sigma^2_{\Delta\hbox{\tiny BMI}})$ are independent.  This assignment
model is complemented by a model describing how BP and $Y$
respond to the prior treatments:
\begin{align}
  \hbox{BP} &= \beta_0+ \beta_1\log(\hbox{Ex}) + \varepsilon_{\hbox{\tiny
  BP}},\label{eq:bp}\\ 
  \hbox{$Y$} &= \delta_0+ \delta_1\log(\hbox{Ex}) + \delta_2\,\Delta \hbox{BMI}
  + \varepsilon_{\hbox{\tiny $Y$}},
\label{eq:y}
\end{align}
where $(\varepsilon_{\hbox{\tiny BP}},\varepsilon_{\hbox{\tiny $Y$}})^t$
are centered bivariate normal and independent of $\varepsilon_{\hbox{\tiny
    Ex}}$ and $\varepsilon_{\Delta \hbox{\tiny BMI}}$.  We denote the
variances of $\varepsilon_{\hbox{\tiny BP}}$ and $\varepsilon_{\hbox{\tiny
    $Y$}}$ by $\sigma_{\hbox{\tiny BP}}^2$ and $\sigma_{\hbox{\tiny
    $Y$}}^2$, respectively, and write $\sigma_{\hbox{\tiny BP,$Y$}}$ for
the possibly non-zero covariance of $\varepsilon_{\hbox{\tiny BP}}$ and
$\varepsilon_{\hbox{\tiny $Y$}}$.  Figure \ref{fig:0robins} shows the path
diagram for this structural equation model.

Equations (\ref{eq:ex}), (\ref{eq:bmi}) and (\ref{eq:bp}) simply specify
conditional expectations that can be estimated in regressions.  However,
this is not the case in general with (\ref{eq:y}).  Instead,
\[
E\left[\hbox{$Y$} \mid \log(\hbox{Ex}), \Delta\hbox{BMI}\right] =
\bar\delta_0+ \bar\delta_1\log(\hbox{Ex}) + \bar\delta_2\,\Delta \hbox{BMI}
\]
with
\begin{align}
  \bar\delta_1 &= \delta_1
  -\frac{\gamma_2\sigma_{\hbox{\tiny
  BP,$Y$}}(\beta_1\gamma_2+\gamma_1)}{\gamma_2^2\sigma_{\hbox{\tiny
  BP}}^2+\sigma_{\hbox{\tiny $\Delta$BMI}}^2},\\ 
  \bar\delta_2 &= \delta_2+
  \frac{\gamma_2\sigma_{\hbox{\tiny BP,$Y$}}}{\gamma_2^2\sigma_{\hbox{\tiny
  BP}}^2+\sigma_{\hbox{\tiny $\Delta$BMI}}^2}. 
\end{align}
We see that $\delta_1$ and $\delta_2$ would have an interpretation as
regression coefficients if: (i) the assignment of $\Delta$BMI did not
depend on BP (i.e., $\gamma_2=0$) and thus both treatments were completely
randomized, or (ii) there were no dependence between
$\varepsilon_{\hbox{\tiny $Y$}}$ and $ \varepsilon_{\hbox{\tiny BP}}$
(i.e., $\sigma_{\hbox{\tiny BP,$Y$}} = 0$).  Similarly, in
$E\left[\hbox{$Y$} \mid \log(\hbox{Ex}), \hbox{BP},
  \Delta\hbox{BMI}\right]$, the coefficient of $\Delta$BMI is equal to
$\delta_2$ but the coefficient for $\log(\hbox{Ex})$ is
$\delta_1-\beta_1\sigma_{\hbox{\tiny BP,$Y$}}/\sigma_{\hbox{\tiny BP}}^2$.

In this paper we consider likelihood-based methods for fitting a large
class of structural equation models that includes the one given by
(\ref{eq:ex})-(\ref{eq:y}) and can be used for consistent estimation of
parameters such as $\delta_1$.  For alternative semi-parametric methods,
see
\cite{robins:1999} and \cite{gill:robins:causal:2001}.

\subsection{Challenges in Structural Equation Modelling}
\label{sec:challenges}

A number of mathematical and statistical problems arise in the normal
linear models associated with systems of structural equations:

\begin{enumerate}
\item \label{item:equiv} Different path diagrams may induce the same
  statistical model, i.e., family of multivariate normal distributions.
  Such {\em model equivalence\/} occurs, for example, for the two path
  diagrams $1\rightarrow 2$ and $1\leftarrow 2$, which differ substantively
  by the direction of the cause-effect relationship.  The two associated
  statistical models, however, are identical, both allowing for correlation
  between the two variables.
\item \label{item:interpret} In many important special cases the path
  coefficient associated with a directed edge $i\rightarrow j$ has a
  population interpretation as a regression coefficient in a regression of
  $j$ on a set of variables including $i$.  However, as seen already in
  \S\ref{sec:motivating-example}, this interpretation is not valid in
  general.
\item \label{item:identify} The parameters of the model may not be
  identifiable, so two different sets of parameter values may lead to the
  same population distribution; for an early review of this problem see
  \cite{fisher:identification:1966}.
\item \label{item:cef} The set of parameterized covariance matrices may
  contain `singularities' at which it cannot be approximated locally by a
  linear space.  At `singular' points, $\chi^2$ and normal approximation to
  the distribution of likelihood ratio tests and maximum likelihood
  estimators (MLE) may not be valid; see for instance \cite{drton:lrt}.
\item \label{item:mle} Iterative procedures are typically required for
  maximization of the likelihood function, which for some models can be
  multimodal \citep{drton:2004}.  Such multimodality typically occurs in
  small samples or under model misspecification.
\end{enumerate}

The problems listed may arise in models without unobserved variables and
become only more acute in latent variable models.  They are challenging in
full generality, but significant progress has been made in special cases
such as recursive linear models with uncorrelated errors, which are also
known as directed acyclic graph (DAG) models or `Bayesian' networks
\citep{lau:bk,pearl:1988}.  A normal DAG model is equivalent to a series of
linear regressions, is always identified and has standard asymptotics.
Under simple sample size conditions, the MLE exists almost surely and is a
rational function of the data.  Graphical modelling theory also solves the
equivalence problem (i) by characterizing all DAGs that induce the same
statistical model \citep{anderson:characterization:1997}.  For more recent
progress on the equivalence problem (\ref{item:equiv}) see
\cite{ali:richardson:spirtes:zhang:towards:2005} and \cite{zhang:spirtes:2005}.

\subsection{New Contribution}
\label{sec:contributions}

The requirement of uncorrelated errors may be overly restrictive in many
settings.  While arbitrary correlation patterns over the errors may yield
rather ill-behaved statistical models, there are subclasses of models with
correlated errors in which some of the nice properties of DAG models are
preserved; compare \cite{mcdonald:what:2002}.  In this paper we consider
path diagrams in which there are no directed cycles and no `double' edges
of the form $i\bidir j$ (compare Def.~\ref{def:acyclic} and
\ref{def:bows}).  Since such double edges have been called `bows', we call
this class {\it bow-free acyclic path diagrams} (BAPs).  An example of a
BAP arose in our motivating example in \S\ref{sec:motivating-example}; see
Figure~\ref{fig:0robins}. While instrumental variable models, which are
much studied in economics, contain bows, most models in other social
sciences are based on BAPs.  For instance, all path diagrams in
\cite{bollen:89} are BAPs.

Bow-free acyclic path diagrams were also considered by \cite{brito:2002}
who showed that the associated normal linear models are almost everywhere
identifiable; see \S\ref{sec:bowfree} for the definition.  This result and
other identification properties of BAP models are reviewed in Section
\ref{sec:graph-view}.
In Section \ref{sec:lik-inf} we give details on likelihood equations and
Fisher-information of normal structural equation models.  This sets the
scene for our main contribution: the {\it Residual Iterative Conditional
  Fitting} (RICF) algorithm for maximization of the likelihood function of
BAP models, which is presented in Section \ref{sec:icf}.  Standard software
for structural equation modelling
currently employs general-purpose optimization routines for this task
\citep{bollen:89}.  Many of these algorithms, however, neglect constraints
of positive definiteness on the covariance matrix and suffer from
convergence problems.
According to \citet{steiger:driving:2001}, failure to converge is `not
uncommon' and presents significant challenges to novice users of existing
software.  In contrast, our RICF algorithm produces positive definite
covariance matrix estimates during all its iterations and has very good
convergence properties, as illustrated in the simulations in
Section~\ref{sec:simulations}.  Further discussion of RICF is provided in
Section \ref{sec:conclusion}.


\section{Normal Linear Models and Path Diagrams}
\label{sec:graph-view}

Let $Y=(Y_i\mid i\in V)\in\RRR^V$ be a random vector, indexed by the finite
set $V$, that follows a multivariate normal distribution $\ND(0,\Sigma)$
with positive definite covariance matrix $\Sigma$.  A zero mean vector is
assumed merely to avoid notational overhead.  The models we consider
subsequently are induced by linear structural equations as follows.

\subsection{Systems of Structural Equations and Path Diagrams}

\begin{figure}[t]
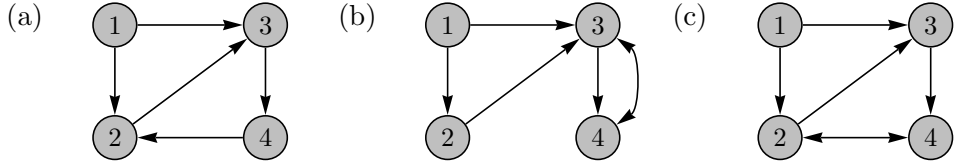

\begin{center}
    \small
    \vspace{1.75cm}\hspace{0.25cm}
  \begin{tabular}{p{4cm}p{4cm}p{4cm}}
    \psset{linewidth=0.6pt,arrowscale=1.5 2,arrowinset=0.1}  
    \psset{fillcolor=lightgray, fillstyle=solid}          
    \newcommand{\myNode}[2]{\circlenode{#1}{\makebox[1.75ex]{#2}}}
    \rput(-0.2,1.6){\normalsize(a)}
    \rput(1, 1.5){\myNode{1}{$1$}}
    \rput(1, 0){\myNode{2}{$2$}} 
    \rput(3, 1.5){\myNode{3}{$3$}} 
    \rput(3, 0){\myNode{4}{$4$}} 
    \ncline{->}{1}{2}
    \ncline{->}{1}{3}
    \ncline{->}{2}{3}
    \ncline{->}{3}{4}
    \ncline{<-}{2}{4}
    & 
    \psset{linewidth=0.6pt,arrowscale=1.5 2,arrowinset=0.1}          
    \psset{fillcolor=lightgray, fillstyle=solid}          
    \newcommand{\myNode}[2]{\circlenode{#1}{\makebox[1.75ex]{#2}}}
    \rput(-0.2,1.6){\normalsize(b)}
    \rput(1, 1.5){\myNode{1}{$1$}}
    \rput(1, 0){\myNode{2}{$2$}} 
    \rput(3, 1.5){\myNode{3}{$3$}} 
    \rput(3, 0){\myNode{4}{$4$}} 
    \psset{fillcolor=lightgray, fillstyle=none}          
    \ncline{->}{1}{2}
    \ncline{->}{1}{3}
    \ncline{->}{2}{3}
    \ncline{->}{3}{4}
    \nccurve[angleA=-35,angleB=35]{<->}{3}{4}
    &
    \psset{linewidth=0.6pt,arrowscale=1.5 2,arrowinset=0.1}          
    \psset{fillcolor=lightgray, fillstyle=solid}          
    \newcommand{\myNode}[2]{\circlenode{#1}{\makebox[1.75ex]{#2}}}
    \rput(-0.2,1.6){\normalsize(c)}
    \rput(1, 1.5){\myNode{1}{$1$}}
    \rput(1, 0){\myNode{2}{$2$}} 
    \rput(3, 1.5){\myNode{3}{$3$}} 
    \rput(3, 0){\myNode{4}{$4$}} 
    \ncline{->}{1}{2}
    \ncline{->}{1}{3}
    \ncline{->}{2}{3}
    \ncline{->}{3}{4}
    \ncline{<->}{2}{4}
  \end{tabular}
\end{center}
\caption{\label{fig:2simple} Path diagrams that
  are (a) cyclic, (b) acyclic but not bow-free, (c) acyclic and bow-free.
  Only path diagram (c) yields a curved exponential
  family.}
\end{figure}

Let $\{\pa(i) \mid i\in V\}$ and $\{\spo(i)\mid i\in V\}$ be two families
of index sets satisfying $i\not\in \pa(i) \cup \spo(i)$ for all $i\in V$.
Moreover, let the family $\{\spo(i) \mid
i\in V\}$ satisfy the symmetry condition that $j\in \spo(i)$ if and only if
$i\in \spo(j)$.  These two families determine a system of structural
equations
\begin{equation}
\label{eq:regeqn}
Y_i=\lsum_{j\in \pa (i)}\beta_{ij}\,Y_j+\veps_i,\quad i\in V,
\end{equation}
whose zero-mean errors $\veps_i$ and $\veps_j$ are uncorrelated if
$i\not\in \spo(j)$, or equivalently, $j\not\in \spo(i)$.  The equations in
(\ref{eq:regeqn}) correspond to a {\em path diagram\/}, that is, a mixed
graph $G$ featuring both {\em directed} ($\to$) and {\em bi-directed}
($\bi$) edges but no edges from a vertex $i$ to itself (see
Figures~\ref{fig:0robins} and \ref{fig:2simple}).  The vertex set of $G$ is
the index set $V$, and $G$ contains the edge $j\to i$ if and only if $j\in
\pa(i)$, and the edge $j\bi i$ if and only if $j\in \spo(i)$ (or
equivalently, $i\in \spo(j)$).  Subsequently, we will exploit the path
diagram representation of \eqref{eq:regeqn}.  If $i\to j$ is an edge in
$G$, then we call $i$ a {\em parent\/} of $j$, and if $i\bi j$ is in $G$
then $i$ is referred to as a {\em spouse\/} of $j$.  Thus $\pa(i)$,
$\spo(i)$ are, respectively, the sets of parents and spouses of $i$.

Let $G$ be a path diagram and define $\mathbf{B}(G)$ to be the
collection of all $V\times V$ matrices $B=(\beta_{ij})$ that satisfy
\begin{equation}
\label{eq:B-constraints}
j\to i\text{ not in }G\;\Longrightarrow\;\beta_{ij}=0, 
\end{equation}
and are such that $I-B$ is invertible.  Let $\mathbf{P}(V)$ be the cone of
positive definite $V\times V$ matrices and $\mathbf{O}(G)\subseteq
\mathbf{P}(V)$ the set of matrices $\Omega=(\omega_{ij})\in \mathbf{P}(V)$
that satisfy
\begin{equation}
\label{eq:O-constraints}
i\not=j \text{ and } j\bi i\text{ not in }G\;\Longrightarrow\;\omega_{ij}=0.
\end{equation}
(Here and in the sequel, a symbol such as $V$ denotes both a finite set and
its cardinality.)  The system (\ref{eq:regeqn}) associated with the path
diagram $G$ can be written compactly as $Y=B\,Y+\veps$.  If we assume that
$B\in\mathbf{B}(G)$ and that the error covariance matrix
$\Var(\veps)=\Omega$ is in $\mathbf{O}(G)$, then (\ref{eq:regeqn}) has a
unique solution $Y$ that is a multivariate normal random vector with
covariance matrix
\begin{equation}
\label{eq:sigmaBOmega}
\Sigma=\Var(Y)=\Phi_G(B,\Omega):=(\idmat-B)^{-1}\Omega(\idmat-B)^{-t}.
\end{equation}
Here, $\idmat$ is the identity matrix and the superscript `$-t$' stands for
transposition and inversion.

The above considerations lead to the following definition of a linear model
associated with a path diagram (or equivalently, a system of structural
equations).

\begin{definition}
  The normal linear model $\mathbf{N}(G)$ associated with a path diagram
  $G$ is the family of multivariate normal distributions $\ND(0,\Sigma)$
  with covariance matrix in the set $\mathbf{P}(G)=\left\{ (\idmat-B)^{-1}
    \Omega (\idmat-B)^{-t}\mid B\in\mathbf{B}(G),\,\Omega\in\mathbf{O}(G)
  \right\}$.  We call the map $\Phi_G:\mathbf{B}(G)\times
  \mathbf{O}(G)\to\mathbf{P}(G)$ defined in (\ref{eq:sigmaBOmega}) the
  parameterization map of $\mathbf{N}(G)$.
\end{definition}

\begin{example}\label{ex:cyclic}\rm
  The path diagram $G$ in Figure \ref{fig:2simple}(a) depicts the
  equation system 
  \begin{align*}
    Y_1&=\veps_1,&
    Y_2&=\beta_{21}Y_1+\beta_{24}Y_4+\veps_2,\\
    Y_3&=\beta_{31}Y_1+\beta_{32}Y_2+\veps_3,&
    Y_4&=\beta_{43}Y_3+\veps_4,
  \end{align*}
  where $\veps_1,\veps_2,\veps_3,\veps_4$ are pairwise uncorrelated, that
  is, the matrices $\Omega\in\mathbf{O}(G)$ are diagonal.  This system
  exhibits a circular covariate-response structure as the path diagram
  contains the directed cycle $2\to 3\to 4\to 2$.  This feedback loop is
  reflected in the fact that $\det(I-B)=1-\beta_{24}\beta_{43}\beta_{32}$
  for $B\in\mathbf{B}(G)$.  Therefore, the path coefficients need to
  satisfy that $\beta_{24}\beta_{43}\beta_{32}\not= 1$ in order to lead to
  a positive definite covariance matrix in $\mathbf{P}(G)$.  This example
  is considered in more detail in \cite{drton:lrt}, where it is shown that
  the parameter space $\mathbf{P}(G)$ has singularities that lead to
  non-standard behavior of likelihood ratio tests.
\end{example}

The models considered in the remainder of this paper will not have any
feedback loops, i.e., they have the following structure.

\begin{definition}
  \label{def:acyclic}
  A path diagram $G$ and its associated normal linear model $\mathbf{N}(G)$
  are {\em recursive\/} or {\em acyclic} if $G$ does not contain directed
  cycles, that is, there do not exist $i,i_1,\dots,i_k\in V$ such that $G$
  features the edges $i\to i_1\to i_2\to\dots\to i_k\to i$.
\end{definition}

We will use the term {\em acyclic\/} rather than {\em recursive\/}, as some
authors have used the term `recursive' for path diagrams that are acyclic
{\it and} contain no bi-directed edges.  Acyclic path diagrams coincide
with the summary graphs of \cite{coxwerm:book}.
If $G$ is acyclic, then the vertices in $V$ can be ordered such that a matrix
$B$ that satisfies (\ref{eq:B-constraints}) is lower-triangular.  It
follows that
\begin{equation}
  \label{lem:phipoly}
  \det(I-B)=1.
\end{equation}
In particular, $\idmat-B$ is invertible for any choice of the path
coefficients $\beta_{ij}$, $j\to i$ in $G$, and the parameterization map
$\Phi_G$ is a polynomial map.

\subsection{Bow-free Acyclic Path Diagrams (BAPs)}
\label{sec:bowfree}

The normal linear model $\mathbf{N}(G)$ associated with a path diagram $G$
is {\em everywhere identifiable} if the parameterization map $\Phi_G$ is
one-to-one, i.e., for all $B_0\in\mathbf{B}(G)$ and
$\Omega_0\in\mathbf{O}(G)$ it holds that
\begin{equation}
  \label{eq:ident}
  \Phi_G(B,\Omega)=\Phi_G(B_0,\Omega_0)\quad\Longrightarrow\quad
  B=B_0\text{ and }\Omega=\Omega_0.
\end{equation}
If there exists a Lebesgue null set $N_G\subseteq
\mathbf{B}(G)\times\mathbf{O}(G)$ such that (\ref{eq:ident}) holds for all
$(B_0,\Omega_0)\not\in N_G$, then we say that $\mathbf{N}(G)$ is {\em
  almost everywhere identifiable}.

Acyclic path diagrams may contain {\em bows}, that is, double edges
$i\bidir j$.  It is easy to see that normal linear models associated with
path diagrams with bows are never everywhere identifiable.  However, they
may sometimes be almost everywhere identifiable as is the case for the next
example.  This example illustrates that almost everywhere identifiability
is not enough to ensure regular behaviour of statistical procedures.

\begin{example}\label{ex:bowed}\rm
  The path diagram $G$ in Figure \ref{fig:2simple}(b) features the bow
  $3\bidir 4$.  The associated normal linear model $\mathbf{N}(G)$ is also
  known as an instrumental variable model.  The 9-dimensional parameter
  space $\mathbf{P}(G)$ is part of the hypersurface defined by the
  vanishing of the so-called {\em tetrad}
  $\sigma_{13}\sigma_{24}-\sigma_{14}\sigma_{23}$.  It follows that the
  model $\mathbf{N}(G)$ lacks regularity because the tetrad hypersurface
  has singularities at points $\Sigma\in\mathbf{P}(G)$ with
  $\sigma_{13}=\sigma_{14}=\sigma_{23}=\sigma_{24}=0$.  These singularities
  occur if and only if $\beta_{31}=\beta_{32}=0$, and correspond to points
  at which the identifiability property in (\ref{eq:ident}) fails to hold.
  This lack of smoothness expresses itself statistically, for example, when
  testing the hypothesis $\beta_{31}=\beta_{32}=0$ in model
  $\mathbf{N}(G)$.  Using the techniques in \cite{drton:lrt}, the
  likelihood ratio statistic for this problem can be shown to have
  non-standard behavior with a large-sample limiting distribution that is
  given by the larger of the two eigenvalues of a $2\times 2$-Wishart
  matrix with 2 degrees of freedom and scale parameter the identity matrix.
\end{example}

\begin{definition}
  \label{def:bows}
  A path diagram $G$ and its associated normal linear model $\mathbf{N}(G)$
  are {\em bow-free\/} if $G$ contains at most one edge between any pair of
  vertices.  If $G$ is bow-free and acyclic, we call it a bow-free acyclic
  path diagram (BAP).
\end{definition}

As stressed in the introduction, BAPs are widespread in applications.
Examples are shown in Figures~\ref{fig:0robins}, \ref{fig:2simple}(c) and
\ref{fig:sur}.  Contrary to some path diagrams with bows, the normal linear
models assciated with BAPs are always at least almost everywhere
identifiable.
\begin{theorem}[\cite{brito:2002}]
  \label{thm:ident}
  If $G$ is a BAP, then the normal linear model $\mathbf{N}(G)$ is almost
  everywhere identifiable.
\end{theorem}
 
Many BAP models are in fact everywhere identifiable.

\begin{theorem}[\cite{richardson:2002}]
  \label{thm:ident2}
  Suppose $G$ is an ancestral BAP, that is, $G$ does not contain an edge
  $i\bi j$ such that there is a directed path $j\to i_1\to\dots \to i_k\to
  i$ that leads from vertex $j$ to vertex $i$.  Then the normal linear
  model $\mathbf{N}(G)$ is everywhere identifiable.
\end{theorem}

The next example shows that the condition in Theorem~\ref{thm:ident2} is
sufficient but not necessary for identification.  The characterization of
the class of BAPs whose associated normal linear models are everywhere
identifiable remains an open problem.

\begin{example}
  \label{ex:fig2c}\rm
  The BAP $G$ in Figure~\ref{fig:2simple}(c) is not ancestral because it
  contains the edges $4\bi 2\to 3\to 4$.  Nevertheless, the associated
  normal linear model $\mathbf{N}(G)$ is everywhere identifiable, which can
  be shown by identifying the parameters in $B$ and $\Omega$ row-by-row
  following the order $1<2<3<4$.  It is noteworthy that the model
  $\mathbf{N}(G)$ in this example is not a Markov model, that is, a generic
  multivariate normal distribution in $\mathbf{N}(G)$ will exhibit no
  conditional independence relations.  Instead, the entries of covariance
  matrices $\Sigma=(\sigma_{ij})\in\mathbf{P}(G)$ satisfy that
  \begin{align}
    \label{eq:fig2c-invariant}
    (\sigma_{11}\sigma_{22}-\sigma_{12}^2)
    (\sigma_{14}\sigma_{33}-\sigma_{13}\sigma_{34})
    =
    (\sigma_{13}\sigma_{24}-\sigma_{14}\sigma_{23} ) 
    (\sigma_{12}\sigma_{13} - \sigma_{11}\sigma_{23}).
  \end{align}
  The constraint in (\ref{eq:fig2c-invariant}) has a nice interpretation.
  Let $(Y_1,\dots,Y_4)$ have (positive definite) covariance matrix
  $\Sigma=(\sigma_{ij})$, and define $e_2=Y_2-\sigma_{21}/\sigma_{11}Y_1$ to
  be the residual in the regression of $Y_2$ on $Y_1$.  Then
  (\ref{eq:fig2c-invariant}) holds for $\Sigma$ if and only if $Y_1$ and
  $Y_4$ are conditionally independent given $e_2$ and $Y_3$.
\end{example}

An inspection of the proof of Theorem~\ref{thm:ident} given in
\cite{brito:2002} reveals the following fact.
\begin{lemma}
  \label{lem:cef}
  If the normal linear model $\mathbf{N}(G)$ associated with a BAP $G$ is
  everywhere identifiable, then the (bijective) parameterization map
  $\Phi_G$ has an inverse that is a rational map with no pole on
  $\mathbf{P}(G)$.
\end{lemma}
  
By~(\ref{lem:phipoly}), the parameterization map $\Phi_G$ for a BAP $G$ is
polynomial and thus smooth.  If $\Phi_G^{-1}$ is rational and without pole,
then the image of $\Phi_G$, i.e., $\mathbf{P}(G)$ is a smooth manifold
\cite[see e.g.][II.4]{edwards:calc}.  This has an important consequence.

\begin{corollary}
  \label{cor:cef}
  If the normal linear model $\mathbf{N}(G)$ associated with a BAP $G$ is
  everywhere identifiable, then $\mathbf{N}(G)$ is a curved exponential
  family.
\end{corollary}

The theory of curved exponential families is discussed in \cite{kassvos}.
It implies in particular that maximum likelihood estimators in curved
exponential families are asymptotically normal, and that likelihood ratio
statistics comparing two such families are asymptotically chi-square
regardless of where in the null hypothesis a true parameter is located.
Unfortunately, however, Lemma~\ref{lem:cef} and Corollary~\ref{cor:cef} do
not hold for every BAP.

\begin{figure}[t]
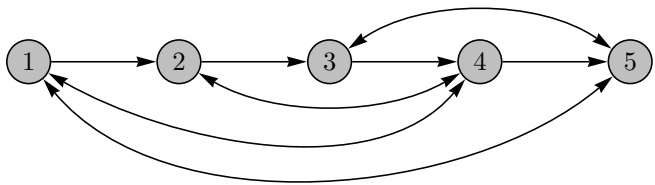

  \centering
    \small
    \vspace{2.5cm}\hspace{-9cm}
    \psset{linewidth=0.6pt,arrowscale=1.5 2,arrowinset=0.1}  
    \psset{fillcolor=lightgray, fillstyle=solid}          
    \newcommand{\myNode}[2]{\circlenode{#1}{\makebox[1.75ex]{#2}}}
    \rput(1, 1.5){\myNode{1}{$1$}}
    \rput(3, 1.5){\myNode{2}{$2$}} 
    \rput(5, 1.5){\myNode{3}{$3$}} 
    \rput(7, 1.5){\myNode{4}{$4$}} 
    \rput(9, 1.5){\myNode{5}{$5$}} 
    \psset{fillcolor=lightgray, fillstyle=none}          
    \ncline{->}{1}{2}
    \ncline{->}{2}{3}
    \ncline{->}{3}{4}
    \ncline{->}{4}{5}
    \nccurve[angleA=-30,angleB=230]{<->}{1}{4}
    \nccurve[angleA=-50,angleB=220]{<->}{1}{5}
    \nccurve[angleA=-30,angleB=210]{<->}{2}{4}
    \nccurve[angleA=35,angleB=145]{<->}{3}{5}
  \caption{Bow-free acyclic path diagram whose associated normal linear
    model is almost, but not everywhere, identifiable.  The model is not a
    curved exponential family.\label{fig:bap-singular}}
\end{figure}

\begin{example}
  \label{ex:bap-singular}\rm
  The normal linear model associated with the BAP $G$ in
  Figure~\ref{fig:bap-singular} is not a curved exponential family.  In
  this model the identifiability property in (\ref{eq:ident}) breaks down
  if and only if $(B,\Omega)$ satisfy that
  \begin{align*}
    \beta_{21}\omega_{14}\omega_{24}-\omega_{2}\omega_{4}+\omega_{24}^2=0, \quad
    \beta_{32}\beta_{43}\omega_{2}+\omega_{24}=0.
  \end{align*}
  It can be shown that the covariance matrices $\Phi_G(B,\Omega)$
  associated with this set of parameters yield points at which the
  13-dimensional set $\mathbf{P}(G)$ has more than 13 linearly independent
  tangent directions.  Hence, $\mathbf{P}(G)$ is singular at these
  covariance matrices.
\end{example}

\section{Likelihood Inference}
\label{sec:lik-inf}

Suppose a sample indexed by a set $N$ is drawn from a multivariate normal
distribution $\ND(0,\Sigma)$ in the linear model $\mathbf{N}(G)$ associated
with a BAP $G=(V,E)$.  We group the observed random vectors as columns in
the $V\times N$ matrix $Y$ such that $Y_{im}$ represents the observation of
variable $i$ on subject $m$.  Having assumed a zero mean vector, we define
the empirical covariance matrix as
\begin{equation}\label{ch5empcov}
  S=(s_{ij})=\frac{1}{N} YY^t\in\RRR^{V\times V}.
\end{equation}
Assuming $N\ge V$, the matrix $S$ is positive definite with probability
one.  Here, the symbols $N$ and $V$ are used to also denote set
cardinalities. Models with unknown mean vector $\mu\in \RRR^V$ can be
treated by estimating $\mu$ by the empirical mean vector and adjusting the
empirical covariance matrix accordingly; $N\geq V+1$ then ensures almost
sure positive definiteness of the empirical covariance matrix.

\subsection{Likelihood Function and Likelihood Equations}
\label{sec:lik-fct}

Given observations $Y$ with empirical covariance matrix $S$, the
log-likelihood function $\ell : \mathbf{B}(G)\times \mathbf{O}(G)\to \RRR$
of the model $\mathbf{N}(G)$ takes the form
\begin{equation}
\label{eq:ellBO}
\ell(B,\Omega)
=-\frac{N}{2}\log\det(\Omega)
-\frac{N}{2}\tr\big[(\idmat-B)^t\Omega^{-1}(\idmat-B)S\big].
\end{equation}
Here we ignored an additive constant and used that $\det(\idmat-B)=1$ if
$B\in\mathbf{B}(G)$; compare (\ref{lem:phipoly}).  Let
$\beta=(\beta_{ij}\mid i\in V,\; j\in\pa(i))$ and $\omega=(\omega_{ij}\mid
i\leq j,\; j\in\spo(i)\text{ or } i=j)$ be the vectors of unconstrained
elements in $B$ and $\Omega$.  Let $P$ and $Q$ be the matrices with entries
in $\{0,1\}$ that satisfy $\vecc(B)=P\beta$ and $\vecc(\Omega)=Q\omega$,
respectively, where $\vecc(A)$ refers to stacking of the columns of the
matrix $A$.  Taking the first derivatives of $\ell(B,\Omega)$ with respect
to $\beta$ and $\omega$ we obtain the likelihood equations.

\begin{proposition}
  \label{prop:likeqns}
  The likelihood equations of the normal linear model $\mathbf{N}(G)$
  associated with a BAP $G$ can be written as
  \begin{align}
    \label{eq:likeqbeta}
    P^t\vecc\big(\Omega^{-1}(\idmat-B)S\big)
    =P^t\vecc\big(\Omega^{-1}S\big)
    -P^t\big(S\otimes\Omega^{-1}\big)P\beta&=0,\\
    \label{eq:likeqomega}
    Q^t\vecc\big(\Omega^{-1}-\Omega^{-1}
    (\idmat-B)S(\idmat-B)^t\Omega^{-1}\big)&=0,
  \end{align}
  where  $\otimes$ denotes the Kronecker product.
\end{proposition}

In general, the likelihood equations need to be solved iteratively. One
possible approach proceeds by alternately solving \eqref{eq:likeqbeta} and
\eqref{eq:likeqomega} for $\beta$ and $\omega$, respectively.  
For fixed $\omega$, \eqref{eq:likeqbeta} is a linear equation in $\beta$
and easily solved.  For fixed $\beta$, \eqref{eq:likeqomega} constitutes
the likelihood equations of a multivariate normal covariance model for
$\varepsilon=(\idmat-B)Y$, which is specified by requiring that
$\Omega_{ij}=0$ whenever the edge $i\bi j$ is not in $G$.  The solution of
\eqref{eq:likeqomega} requires, in general, another iterative method.  As
an alternative to this nesting of two iterative methods, we propose in
Section \ref{sec:icf} a method that solves \eqref{eq:likeqbeta} and
\eqref{eq:likeqomega} in joint updates of $\beta$ and $\omega$.

\subsection{Fisher-Information}
\label{sec:fish-info}

Large-sample confidence intervals for $(\beta,\omega)$ can be obtained by
approximating the distribution of the MLE $(\hat\beta,\hat\omega)$ by the
normal distribution with mean vector $(\beta,\omega)$ and covariance matrix
$\frac{1}{N}\,I(\beta,\omega)^{-1}$.  Here, $I(\beta,\omega)$ denotes the
Fisher-information, which, as shown in
Appendix~\ref{sec:app-proof-identify}, is of the following form.

\begin{proposition} \label{prop:fisherinfo}
  The (expected) Fisher-information of the normal linear model
  $\mathbf{N}(G)$ associated with a BAP $G$ is
  \[
  I(\beta,\omega)=
  \begin{pmatrix}
    P^t\big(\Sigma\otimes\Omega^{-1}\big)P&
    P^t\big[(\idmat-B)^{-1}\otimes\Omega^{-1}\big]Q\\
    Q^t\big[(\idmat-B)^{-t}\otimes\Omega^{-1}\big]P&
    \frac{1}{2}\,Q^t\big(\Omega^{-1}\otimes\Omega^{-1}\big)Q
  \end{pmatrix}.
  \]
\end{proposition}

The Fisher-information in Proposition \ref{prop:fisherinfo} need not be
block-diagonal, in which case the estimation of the covariances $\omega$
affects the asymptotic variance of the MLE $\hat\beta$.  However, this does
not happen for {\em bi-directed chain graphs\/}, which form one of the
model classes discussed in \cite{wermuth:2004}.  A path diagram $G$ is a
bi-directed chain graph if its vertex set $V$ can be partitioned into
disjoint subsets $\tau_1,\dots,\tau_T$, known as {\it chain components},
such that all edges in each subgraph $G_{\tau_t}$ are bi-directed and edges
between two subsets $\tau_s\neq\tau_t$ are directed, pointing from $\tau_s$
to $\tau_t$, if $s<t$.  Since bi-directed chain graphs are ancestral graphs
the associated normal linear models are everywhere identifiable.

\begin{proposition}
  \label{prop:chaingraph}
  For a BAP $G$, the following two statements are equivalent:
  \begin{enumerate}
  \item[(i)] For all underlying covariance matrices
    $\Sigma\in\mathbf{P}(G)$, the MLEs of the parameter vectors $\beta$ and
    $\omega$ of the normal linear model $\mathbf{N}(G)$ are asymptotically
    independent.
  \item[(ii)] The path diagram $G$ is a bi-directed chain graph.
  \end{enumerate}
\end{proposition}

A proof of Proposition \ref{prop:chaingraph} is given in
Appendix~\ref{sec:app-proof-identify}.  This result is an instance of the
asymptotic independence of mean and natural parameters in mixed
parameterization of exponential families \citep{barndorff-nielsen:1978}.

\section{Residual Iterative Conditional Fitting}
\label{sec:icf}

We now present an algorithm for computing the MLE in the normal linear
model $\mathbf{N}(G)$ associated with a BAP.  The algorithm extends the
{\em iterative conditional fitting} (ICF) procedure of \citet{icfsanjay},
which is for path diagrams with exclusively bi-directed edges.

Let $Y_i\in\RRR^N$ denote the $i$-th row of the observation matrix $Y$ and
$Y_{-i}= Y_{V\setminus\{i\}}$ the $(V\setminus\{i\})\times N$ submatrix of
$Y$.  The ICF algorithm proceeds by repeatedly iterating through all
vertices $i\in V$ and carrying out three steps: (i) fix the marginal
distribution of $Y_{-i}$, (ii) fit the conditional distribution of $Y_i$
given $Y_{-i}$ under the constraints implied by the model $\mathbf{N}(G)$,
and (iii) obtain a new estimate of $\Sigma$ by combining the estimated
conditional distribution $(Y_i\mid Y_{-i})$ with the fixed marginal
distribution of $Y_{-i}$.
The crucial point is then that for path diagrams containing only
bi-directed edges, the problem of fitting the conditional distribution for
$(Y_i\mid Y_{-i})$ under the constraints of the model can be rephrased as a
least squares regression problem. Unfortunately, the consideration of the
conditional distribution of $(Y_i\mid Y_{-i})$ is complicated for path
diagrams that contain also directed edges.  However, as we show below, the
directed edges can be `removed' by consideration of {\em residuals\/},
which here refers to estimates of the error terms $\veps=(\idmat-B)Y$.
Since it is based on this idea, we give our new extended algorithm the name
{\em Residual Iterative Conditional Fitting\/} (RICF).

\subsection{The RICF Algorithm}
\label{sec:general-algorithm}

The main building block of the new algorithm is the following decomposition
of the log-likelihood function.  We adopt the shorthand notation $X_C$ for
the $C\times N$ submatrix of a $D\times N$ matrix $X$, where $C\subseteq
D$.

\begin{theorem}
  \label{thm:ricfdecomp}
  Let $G$ be a BAP and $i\in V$.  Let $\|x\|^2=x^tx$ and define
  \begin{equation}
    \label{eq:condvar}
    \omega_{ii.-i}=
    \omega_{ii}-\Omega_{i,-i}\Omega_{-i,-i}^{-1}\Omega_{-i,i}
  \end{equation}
  to be the conditional variance of $\veps_i$ given $\veps_{-i}$; recall
  that $\Omega_{-i,-i}^{-1}=(\Omega_{-i,-i})^{-1}$.  Then the
  log-likelihood function $\ell(B,\Omega)$ of the model $\mathbf{N}(G)$ can
  be decomposed as
  \begin{align*}
    \ell(B,\Omega)
    =&-\frac{N}{2}\log\omega_{ii.-i}
    -\frac{1}{2\omega_{ii.-i}}\big\|Y_{i}-B_{i,\pa(i)}Y_{\pa(i)}
    -\Omega_{i,\spo(i)}
    (\Omega_{-i,-i}^{-1}\,\varepsilon_{-i})_{\spo(i)}\big\|^2\\
    &\quad\quad-\frac{N}{2}\log\det(\Omega_{-i,-i})
    -\frac{1}{2}\tr\big(\Omega_{-i,-i}^{-1}\veps_{-i}\veps_{-i}^t\big).
  \end{align*}
\end{theorem}
\begin{proof}
  Forming $\veps=(\idmat-B)\,Y$, we rewrite
  (\ref{eq:ellBO}) as
  \begin{equation}
    \label{icf:likveps}
    \ell(B,\Omega)
    =-\frac{N}{2}\log\det(\Omega)-\frac{1}{2}\tr\big(\Omega^{-1}\veps\veps^t\big)
    =:\ell(\Omega\mid \veps).
  \end{equation}
  Using the inverse variance lemma \citep[Prop.~5.7.3]{whittaker:bk},
  we partition $\Omega^{-1}$ as
  \[
  \begin{split}
    &\begin{pmatrix}
      \omega_{ii}&\Omega_{i,-i}\\
      \Omega_{-i,i}&\Omega_{-i,-i}
    \end{pmatrix}^{-1}\!\!\!
    =\begin{pmatrix}
      \omega_{ii.-i}^{-1}&
      -\omega_{ii.-i}^{-1}\Omega_{i,-i}\Omega_{-i,-i}^{-1}\\
      -\Omega_{-i,-i}^{-1}\Omega_{-i,i}\omega_{ii.-i}^{-1}& 
      \Omega_{-i,-i}^{-1}+
      \Omega_{-i,-i}^{-1}
      \Omega_{-i,i}\omega_{ii.-i}^{-1}\Omega_{i,-i}\Omega_{-i,-i}^{-1} 
    \end{pmatrix}.
  \end{split}
  \]
  We obtain that the
  log-likelihood function in \eqref{icf:likveps} equals
  \begin{align}
    \nonumber \ell(\Omega\mid \veps)\;
    =&-\frac{N}{2}\log\omega_{ii.-i}-\frac{1}{2\omega_{ii.-i}}\big\|
    \veps_{i}-\Omega_{i,-i}\Omega_{-i,-i}^{-1}\veps_{-i}\big\|^2\\
    \label{eq:RHS}
    &\quad\quad-\frac{N}{2}\log\det(\Omega_{-i,-i})
    -\frac{1}{2}\tr\big(\Omega_{-i,-i}^{-1}\veps_{-i}\veps_{-i}^t\big).
  \end{align}
  By definition, $\veps_i=Y_i-B_{i,\pa(i)}Y_{\pa(i)}$.  Moreover, under the
  restrictions (\ref{eq:O-constraints}), 
  \[
  \Omega_{i,-i}\Omega_{-i,-i}^{-1}\veps_{-i}=
  \Omega_{i,\spo(i)}
    (\Omega_{-i,-i}^{-1}\,\varepsilon_{-i})_{\spo(i)},
  \]
  which yields the claimed decomposition.
\end{proof}
  
The log-likelihood decomposition (\ref{icf:likveps}) is essentially based
on the decomposition of the joint distribution of $\veps$ into the marginal
distribution of $\veps_{-i}$ and the conditional distribution
$(\veps_{i}\mid \veps_{-i})$.  This leads to the idea of building an
iterative algorithm whose steps are based on fixing the marginal
distribution of $\veps_{-i}$ and estimating a conditional distribution.  In
order to fix the marginal distribution $\veps_{-i}$ we need to fix the
submatrix $\Omega_{-i,-i}$ comprising all but the $i$-th row and column of
$\Omega$ as well as the submatrix $B_{-i,V}$, which comprises all but the
$i$-th row of $B$.
With $\Omega_{-i,-i}$ and $B_{-i,V}$ fixed
we can compute $\veps_{-i}$ as well as the {\em pseudo-variables}, defined by
\begin{equation}
\label{icf:pseudo}
Z_{-i} = \Omega_{-i,-i}^{-1}\,\varepsilon_{-i}.
\end{equation}
From (\ref{icf:likveps}), it now becomes apparent that, for fixed
$\Omega_{-i,-i}$ and $B_{-i,V}$, the maximization of the
log-likelihood function $\ell(B,\Omega)$ can be solved by maximizing
the function
\begin{multline}
  \big((\beta_{ij})_{j\in\pa(i)},(\omega_{ik})_{k\in\spo(i)},
  \omega_{ii.-i})\mapsto\\
  -\frac{N}{2}\log\omega_{ii.-i}
  -\frac{1}{2\omega_{ii.-i}}\big\|Y_{i}-\sum_{j\in\pa(i)}\beta_{ij}Y_j
  -\sum_{k\in\spo(i)}\omega_{ik}Z_{k}\big\|^2
\label{eq:1}
\end{multline}
over $\RRR^{\pa(i)}\times\RRR^{\spo(i)}\times (0,\infty)$.  The maximizers
of (\ref{eq:1}), however, are the least squares estimates in the regression
of
$Y_i$ on both $(Y_j\mid j\in\pa(i))$ and $(Z_k\mid k\in\spo(i))$.

Employing the above observations, the {\em RICF algorithm} for computing
the MLE $(\hat B,\hat\Omega)$ repeats the following steps for each $i\in
V$: \renewcommand{\theenumi}{\arabic{enumi}}
\renewcommand{\labelenumi}{\theenumi.}
\begin{enumerate}
\item \label{item:1} Fix $\Omega_{-i,-i}$ and $B_{-i,V}$, and compute
  residuals $\varepsilon_{-i}$ and pseudo-variables $Z_{\spo(i)}$;
\item\label{item:step4} Obtain least squares estimates of $\beta_{ij}$,
  $j\in\pa(i)$, $ \omega_{ik}$, $k\in\spo(i)$, and $ \omega_{ii.-i}$ by
  regressing response variable $Y_i$ on the covariates $Y_j$, $j\in\pa(i)$ and
  $Z_k$, $k\in\spo(i)$;
\item \label{item:2} Compute an estimate of
  $\omega_{ii}=\omega_{ii.-i}+\Omega_{i,-i}\Omega_{-i,-i}^{-1}\Omega_{-i,i}$
  using the new estimates and the fixed parameters;
  compare (\ref{eq:condvar}).
\end{enumerate}
After steps \ref{item:1} to \ref{item:2}, we move on to the next vertex in
$V$.  After the last vertex in $V$ we return to consider the first vertex.
The procedure is continued until convergence.

\begin{figure}[t]
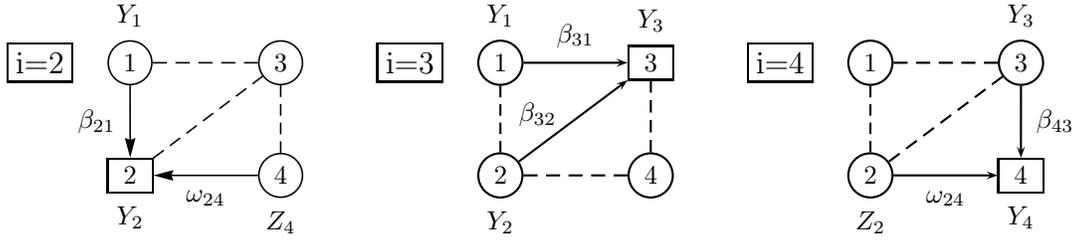

        \begin{center}
        \begin{minipage}{7cm}
          \small \vspace{2cm} \hspace{-3.25cm}
          \newcommand{\myNode}[2]{\circlenode{#1}{\makebox[1.75ex]{#2}}}
  \begin{tabular}{p{4.5cm}p{4.5cm}p{4.5cm}}
    \psset{linewidth=0.6pt,arrowscale=1.5 2,arrowinset=0.1}          
    \rput(-0.2, 1.5){\large$\fbox{i=2}$}
    \rput(1, 1.5){\myNode{1}{$1$}}
    \rput(1, 0){\rnode{2}{\framebox{\makebox[2.5ex]{$2$}}}}
    \rput(3, 1.5){\myNode{3}{$3$}} 
    \rput(3, 0){\myNode{4}{$4$}} 
    \nput{90}{1}{$Y_1$}
    \nput{270}{2}{$Y_2$}
    \nput{270}{4}{$Z_4$}
    \ncline{->}{1}{2}\Bput{$\beta_{21}$}
    \ncline[linestyle=dashed]{-}{2}{3}
    \ncline[linestyle=dashed]{-}{3}{4}
    \ncline[linestyle=dashed]{-}{1}{3}
    \ncline{<-}{2}{4}\Bput{$\omega_{24}$}
    &
    \rput(-0.2, 1.5){\large$\fbox{i=3}$}
    \rput(1, 1.5){\myNode{1}{$1$}}
    \rput(1, 0){\myNode{2}{$2$}} 
    \rput(3, 1.5){\rnode{3}{\framebox{\makebox[2.5ex]{$3$}}}}
    \rput(3, 0){\myNode{4}{$4$}} 
    \nput{90}{1}{$Y_1$}
    \nput{270}{2}{$Y_2$}
    \nput{90}{3}{$Y_3$}
    \ncline{->}{1}{3}\Aput{$\beta_{31}$}
    \ncline[linestyle=dashed]{-}{1}{2}
    \ncline{->}{2}{3}\naput[npos=0.35,labelsep=0.05]{$\beta_{32}$}
    \ncline[linestyle=dashed]{-}{3}{4}
    \ncline[linestyle=dashed]{-}{2}{4}
    &
    \rput(-0.2, 1.5){\large$\fbox{i=4}$}
    \rput(1, 1.5){\myNode{1}{$1$}}
    \rput(1, 0){\myNode{2}{$2$}} 
    \rput(3, 1.5){\myNode{3}{$3$}} 
    \rput(3, 0){\rnode{4}{\framebox{\makebox[2.5ex]{$4$}}}}
    \nput{270}{2}{$Z_2$}
    \nput{90}{3}{$Y_3$}
    \nput{270}{4}{$Y_4$}
    \ncline[linestyle=dashed]{-}{1}{3}
    \ncline[linestyle=dashed]{-}{1}{2}
    \ncline[linestyle=dashed]{-}{2}{3}
    \ncline{->}{3}{4}\Aput{$\beta_{43}$}
    \ncline{->}{2}{4}\Bput{$\omega_{24}$}
    \end{tabular}
    \vspace{0.25cm}
        \end{minipage}
        \end{center}
\caption{\label{fig:icfexample}
  Illustration of the RICF update steps in Example~\ref{ex:ricf}. The
  structure of each least squares regression is indicated by directed edges
  pointing from the predictor variables to the response variable depicted
  by a square node.  (See text for details.)}
\end{figure}

\begin{example}\label{ex:ricf}\rm
  For illustration of the regressions performed in RICF, we consider the
  normal linear model associated with the BAP $G$ in Figure
  \ref{fig:2simple}(c).  The parameters to be estimated in this model are
  $\beta_{21}$, $\beta_{31}$, $\beta_{32}$, $\beta_{43}$ and $\omega_{11}$,
  $\omega_{22}$, $\omega_{33}$, $\omega_{44}$, $\omega_{24}$.
  
  Vertex 1 in Figure \ref{fig:2simple}(c) has no parents or spouses, and
  its RICF update step consists of a trivial regression.  In other words,
  the variance $\omega_{11}$ is the unconditional variance of $Y_1$ with
  MLE $\hat\omega_{11}=s_{11}$.  For the remaining vertices, the
  corresponding RICF update steps are illustrated in Figure
  \ref{fig:icfexample}, where the response variable $Y_i$ in the $i$-th
  update step is shown as a square node while the remaining variables are
  depicted as circles.  Directed edges indicate variables acting as
  covariates in the least squares regression.  These covariates are
  labelled according to whether the random variable $Y_j$, or the
  pseudo-variable $Z_j$ defined in (\ref{icf:pseudo}), is used in the
  regression.  Note that since $\spo(3)=\emptyset$, repetition of steps
  \ref{item:1}--\ref{item:2} in \S\ref{sec:general-algorithm} is required
  only for $i\in\{2,4\}$.
\end{example}

\renewcommand{\theenumi}{\roman{enumi}}
\renewcommand{\labelenumi}{(\theenumi)}
 
In RICF, the log-likelihood function $\ell(B,\Omega)$ from \eqref{eq:ellBO}
is repeatedly maximized over sections in the parameter space defined by
fixing the parameters $\Omega_{-i,-i}$, and $B_{-i,V}$.  RICF thus is an
iterative partial maximization algorithm and has the following convergence
properties. 

\begin{theorem}
  \label{prop:converge}
  If $G$ is a BAP and the empirical covariance matrix $S$ is positive
  definite, then the following holds:
  \begin{enumerate}
  \item For any starting value $(\hat
    B^0,\hat\Omega^0)\in\mathbf{B}(G)\times \mathbf{O}(G)$, RICF constructs
    a sequence of estimates $(\hat B^s,\hat\Omega^s)_s$ in
    $\mathbf{B}(G)\times \mathbf{O}(G)$ whose accumulation points are local
    maxima or saddle points of the log-likelihood function
    $\ell(B,\Omega)$.  Moreover, evaluating the log-likelihood function at
    different accumulation points yields the same value.
  \item If the normal linear model $\mathbf{N}(G)$ is everywhere
    identifiable and the likelihood equations have only finitely many
    solutions then the sequence $(\hat B^s,\hat\Omega^s)_s$ converges to
    one of these solutions.
\end{enumerate}
\end{theorem}
\begin{proof}
  Let $\ell(\Sigma)$ be the log-likelihood function for the model of all
  centered multivariate normal distributions on $\RRR^V$.  If $S$ is
  positive definite then the set $C$ that comprises all positive definite
  matrices $\Sigma\in\RRR^{V\times V}$ at which $\ell(\Sigma)\ge\ell(\hat
  B^0,\hat\Omega^0)$ is compact.  In particular, the log-likelihood
  function in (\ref{eq:ellBO}) is bounded, and claim (i) can be derived
  from general results about iterative partial maximization algorithms; see
  e.g.~\cite{drtoneichler:2004}.  For claim (ii) note that if
  $\mathbf{N}(G)$ is everywhere identifiable, then the compact set $C$ has
  compact preimage $\phi_G^{-1}(C)$ under the model parameterization map;
  recall Lemma~\ref{lem:cef}.
\end{proof}

\begin{remark}
  \label{rem:diverge}\rm
  If the normal linear model $\mathbf{N}(G)$ associated with a BAP $G$ is
  not everywhere identifiable, then it is possible that a sequence of
  estimates $(\hat B^s,\hat\Omega^s)_s$ produced by RICF diverges and does
  not have any accumulation points.  In these cases, however, the
  corresponding sequence of covariance matrices $\Phi_G(\hat
  B^s,\hat\Omega^s)_s$ still has at least one accumulation point because it
  ranges in the compact set $C$ exhibited in the proof of
  Theorem~\ref{prop:converge}.  Divergence of $(\hat B^s,\hat\Omega^s)_s$
  occurs in two instances in the simulations in \S\ref{sec:simulations};
  compare Table~\ref{tab:sims-gene-expression}.  In both cases, the
  sequence $\Phi_G(\hat B^s,\hat\Omega^s)_s$ converges to a positive
  definite covariance matrix.
\end{remark}

\subsection{Computational Savings in RICF}
\label{sec:savericf}

If $G$ is a DAG, i.e., an acyclic path diagram without bi-directed edges,
then the MLE $(\hat B,\hat\Omega)$ in $\mathbf{N}(G)$ can be found in a
finite number of regressions \citep[e.g.][]{wermuth:1980}.  However, we can
also run RICF.  Since in a DAG, $\spo(i)=\emptyset$ for all $i\in V$, step
\ref{item:step4} of RICF regresses variable $Y_i$ solely on its parents
$Y_j$, $j\in\pa(i)$.  Not involving pseudo-variables that could change from
one iteration to the other, this regression remains the same throughout
different iterations, and RICF converges in one step.

Similarly, for a general BAP $G$, if vertex $i\in V$ has no spouses,
$\spo(i)=\emptyset$, then the MLE of $B_{i,\pa(i)}$ and $\omega_{ii}$ can
be determined by a single iteration of the algorithm.  In other words, RICF
reveals these parameters as being estimable in closed form, namely as
rational functions of the data.  (This occurred for vertex $i=3$ in
Example~\ref{ex:ricf}.)  It follows that, to estimate the remaining
parameters, the iterations need only be continued over vertices $j$ with
$\spo(j)\neq\emptyset$.

For further computational savings note that
$\Omega_{\mathrm{dis}(i),V\setminus (\mathrm{dis}(i)\cup \{i\})}=0$, where
$\mathrm{dis}(i)=\{j \mid j\leftrightarrow \cdots \leftrightarrow i, j\neq
i\}$.  Hence, since $\spo (i)\subseteq \mathrm{dis}(i)$,
\[
(\Omega_{-i,-i}^{-1}\varepsilon_{-i})_{\spo (i)} =
    (\Omega_{\mathrm{dis}(i),\mathrm{dis}(i)}^{-1}\varepsilon_{\mathrm{dis}(i)})_{\spo
    (i)}; 
\]
see \citet[Lemma 3.1.6]{koster:fields} and \citet[Lemma
8.10]{richardson:2002}.  Since $\varepsilon_{\mathrm{dis}(i)} =
Y_{\mathrm{dis}(i)} - B_{\mathrm{dis}(i), \pa(\mathrm{dis}(i))}
Y_{\pa(\mathrm{dis}(i))}$, it follows
that in the RICF update step for vertex $i$ attention can be restricted to
the variables in $\{i\}\cup \pa (i) \cup \mathrm{dis}(i)\cup
\pa(\mathrm{dis}(i))$.

Finally, note that while the RICF algorithm is described in terms of the
entire data matrix $Y$, the least squares estimates computed in its
iterations are clearly functions of the empirical covariance matrix, which
is a sufficient statistic.

\section{Simulation Studies}
\label{sec:simulations}

In order to evaluate the performance of the RICF algorithm we consider two
scenarios.  First, we fit linear models based on randomly generated BAPs to
gene expression data.  This scenario is relevant for model selection tasks,
and we compare RICF's performance in this problem to that of algorithms
implemented in software for structural equation modelling.  Second, we
study how RICF behaves when it is used to fit larger models to data
simulated from the respective model.  In contrast to the first scenario,
the second scenario involves models that generally fit the considered data
well.

\subsection{Gene Expression Data}
\label{subsec:gene-ex-data}

We consider data on gene expression in {\it Arabidopsis thaliana\/} from
\cite{willeetal:2004}.  We restrict attention to 13 genes that belong to
one pathway: DXPS1-3, DXR, MCT, CMK, MECPS, HDS, HDR, IPPI1, GPPS, PPDS1-2.
Data from $n=118$ microarray experiments are available.  We fit randomly
generated BAP models to these data using RICF and two alternative methods.

The BAP models are generated as follows.  For each of the 78 possible pairs
of vertices $i<j$ in $V=\{1,\dots,13\}$ we draw from a multinomial
distribution to generate a possible edge.  The probability for drawing the
edge $i\to j$ is $d$, and the probability for drawing $i\bi j$ is $b$ so
that with probability $1-d-b$ there is no edge between $i$ and $j$.  We
then apply a random permutation to the vertices to obtain the final BAP.
For each of twelve combinations $(d,b)$ with $d=0.05,0.1,0.2,0.3$ and
$b=0.05,0.1,0.2$, we simulate 1000 BAPs.  The expected number of edges thus
varies between 7.8 and 39.

For fitting the simulated BAPs to the gene expression data, we implemented
RICF in the statistical programming environment R \citep{R}.  As
alternatives, we consider the R package `sem' \citep{fox:2006} and the very
widely used software LISREL \citep{lisrel} in its student version 8.7 for
Linux (student versions are free but limited to 15 variables).  Both these
latter programs employ general purpose optimizers, e.g., `sem' makes a call
to the R function `nlm'.

\begin{table}[t]
  \centering
  \caption{Fitting simulated BAPs to gene expression data using RICF, LISREL
    and `sem'.  Each row is based on 1000 simulations.  Running time is 
  average CPU time (in sec.) for the cases in which all three algorithms
  converged. (See text for details.)}
  \medskip
  \label{tab:sims-gene-expression}
  \begin{tabular}{llccrrccccrrr}
    \hline\hline
    &&&
    \multicolumn{3}{c}{No convergence} && \multicolumn{1}{c}{All} 
    & \multicolumn{1}{c}{All} &&
    \multicolumn{3}{c}{Running time}\\ 
    \multicolumn{1}{c}{$d$} & 
    \multicolumn{1}{c}{$b$} && 
    \small RICF &  \small LIS &  \small SEM && 
    \multicolumn{1}{c}{converge} & 
    \multicolumn{1}{c}{agree} && \small
    RICF & \small LIS & \small SEM\\ 
    \cline{1-2}
    \cline{4-6}
    \cline{8-9}
    \cline{11-13}                                                                           
$0.05$ &$0.05$ && 0 &  36 &  47 && 941 & 940 && 0.03 & 0.02 & 1.15 \\
       &$0.1$  && 0 & 177 & 221 && 746 & 739 && 0.09 & 0.03 & 1.58 \\ 
       &$0.2$  && 0 & 499 & 599 && 347 & 333 && 0.21 & 0.04 & 2.71 \\
$0.1$  &$0.05$ && 0 &  32 &  36 && 951 & 949 && 0.04 & 0.03 & 1.58 \\
       &$0.1$  && 0 & 137 & 193 && 786 & 780 && 0.09 & 0.03 & 2.09 \\
       &$0.2$  && 0 & 440 & 610 && 364 & 354 && 0.25 & 0.04 & 3.43 \\
$0.2$  &$0.05$ && 0 &  19 &  39 && 958 & 954 && 0.05 & 0.03 & 2.67 \\
       &$0.1$  && 0 &  91 & 176 && 815 & 808 && 0.13 & 0.04 & 3.34 \\
       &$0.2$  && 1 & 326 & 520 && 461 & 452 && 0.33 & 0.05 & 5.03 \\
$0.3$  &$0.05$ && 0 &  16 &  38 && 960 & 957 && 0.06 & 0.04 & 4.04 \\
       &$0.1$  && 0 &  59 & 136 && 859 & 850 && 0.17 & 0.04 & 4.96 \\
       &$0.2$  && 1 & 225 & 471 && 519 & 490 && 0.40 & 0.06 & 6.97 \\
    \hline\hline
  \end{tabular}
\end{table}



Our simulation results are summarized in
Table~\ref{tab:sims-gene-expression}.  Each row in the table corresponds to
a choice of the edge probabilities $d$ and $b$.  The first three columns
count how often, in 1000 simulations, the three considered fitting routines
failed to converge.  The starting values of LISREL and `sem' were set
according to program defaults, and RICF was started by setting $\hat
B^{(0)}$ and $\hat\Omega^{(0)}$ equal to the MLE in the DAG model
associated with the DAG obtained by removing all bi-directed edges from the
considered BAP.

LISREL and `sem' failed to converge for a rather large number of models.
The LISREL output explained why convergence failed, and virtually all
failures were due to the optimizer converging to matrices that were not
positive definite.  The remedy would be to try new starting values but
doing this successfully in an automated fashion is a challenging problem in
itself.  For RICF convergence failure arose in only two cases.  In both
cases the RICF estimates $(\hat B,\hat\Omega)$ had some diverging entries.
Despite the divergence in $(B,\Omega)$-space, the sequence of associated
covariance matrices $\Phi_G(\hat B,\hat\Omega)$ computed by RICF converged,
albeit very slowly.  Recall that this phenomenon is possible in models that
are almost, but not everywhere, identifiable (Remark~\ref{rem:diverge}).
In these examples LISREL produced similarly divergent sequences with
approximately the same likelihood, and `sem' reported convergence in one
case but gave an estimate whose likelihood was nearly 40 points smaller
than the intermediate estimates computed by LISREL and RICF.

The columns labelled `All converge' and `All agree' in
Table~\ref{tab:sims-gene-expression} show how often all methods converged,
and when this occurred, how often the three computed maxima of the
log-likelihood function were the same up to one decimal place.  Since all
methods are for local maximization, substantial disagreements in the
computed maxima can occur if the likelihood function is multimodal.

Finally, the last three columns give average CPU time use for the cases in
which all three algorithms converge.  These are quoted to show that RICF is
competitive in terms of computational efficiency, but for the following
reasons the precise times should not be used for a formal comparison.  On
the one hand, LISREL is fast because it is compiled code.  This is not the
case for the R-based `sem' and RICF.  On the other hand, the fitting
routines in LISREL and `sem' not only compute the MLE but also produce
various other derived quantities of interest.  This is in contrast to our
RICF routine, which only computes the MLE.

\subsection{Simulated Data}
\label{subsec:simulated-data}

In order to demonstrate how RICF behaves when fitting larger models we use
the algorithm on simulated data.  We consider different choices for the
number of variables $p$ and generate random BAPs according to the procedure
used in \S\ref{subsec:gene-ex-data}.  We limit ourselves to two different
settings for the expected number of edges, choosing $d=0.1$ or $d=0.2$ and
setting $b=d/2$ in each case.  For each BAP $G$, we simulate a covariance
matrix $\Sigma=(\idmat-B)^{-1}\Omega(\idmat-B)^{-t}\in\mathbf{P}(G)$ as
follows.  The free entries in $B\in\mathbf{B}(G)$ and the free off-diagonal
entries in $\Omega\in\mathbf{O}(G)$ are drawn from a $\mathcal{N}(0,1)$
distribution.  The diagonal entries $\omega_{ii}$ are obtained by adding a
draw from a $\chi^2_1$-distribution to the sum of the absolute values of
the off-diagonal entries in the $i$-th row of $\Omega$.  This makes
$\Omega$ diagonally dominant and thus positive definite.  Finally, we draw
a sample of size $n$ from the resulting multivariate normal distribution
$\mathbf{N}(G)$.  For each distribution two cases, namely $n=3p/2$ and
$n=10 p$, are considered to illustrate sample size effects.  For each
combination of $p$, $d$ and $n$, we simulate 500 BAPs and associated data
sets.

\begin{figure}[t]
  \centering
  \hspace{-0.5cm}
  \includegraphics[width=14.5cm]{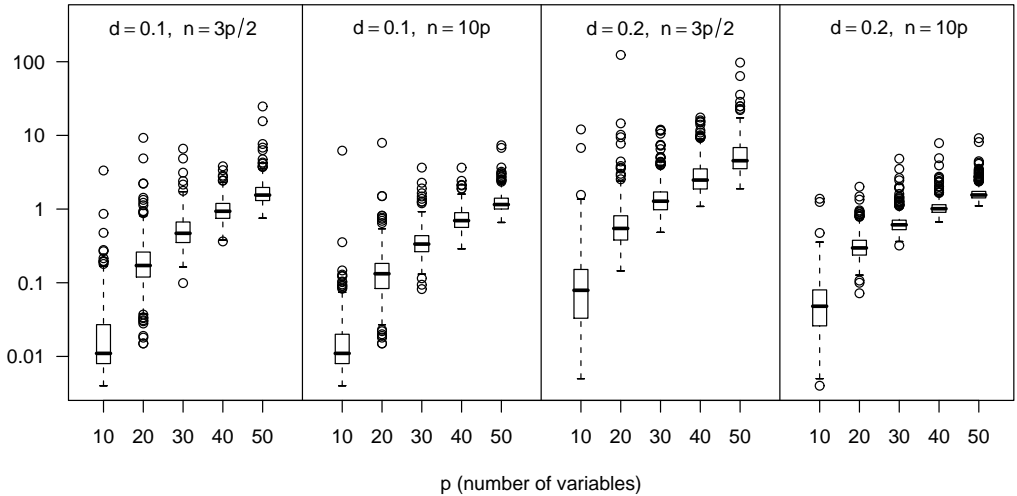}
  \caption{Boxplots of CPU times (in sec.~on $\log_{10}$-scale) used by 
    RICF when fitting BAP models to simulated data.  Each boxplot
    summarizes 500 simulations.  The number of variables is
    denoted by $p$, the sample size is $n$, and the parameter $d$
    determines the expected number of edges of the simulated BAPs (see text
    for details).\label{fig:sim-model}}
\end{figure}

Figure~\ref{fig:sim-model} summarizes the results of our simulations in
boxplots.  As could be expected, the running time for RICF increases with
the number of variables $p$ and the expected number of edges in the BAP
(determined by $d$).  Moreover, the running time decreases for increased
sample size $n$, which is plausible because the empirical covariance matrix
of a larger sample will tend to be closer to the underlying parameter space
$\mathbf{P}(G)$.  The boxplots show that even with $p=50$ variables the
majority of the RICF computations terminate within a few seconds.  However,
there are also a number of computations in which the running time is
considerably longer, though still under two minutes. This occurs in
particular for the denser case with smaller sample size ($d=0.2$ and
$n=3p/2$).

\section{Discussion}
\label{sec:conclusion}

As mentioned in the introduction, normal linear models associated with path
diagrams are employed in many applied disciplines.  The models, also known
as structural equation models, have a long tradition but remain of current
interest in particular due to the more recent developments in causal
inference; compare e.g.~\cite{pearl:causalbk,spirtes:2000}.  Despite their
long tradition, however, many mathematical, statistical and computational
problems about these models remain open.  

The new contribution of this paper is the Residual Iterative Conditional
Fitting (RICF) algorithm for maximum likelihood estimation in BAP models.
Software for computation of MLEs in structural equation models often
employs optimization methods that are not designed to deal with positive
definiteness constraints on covariance matrices.  This can be seen in
particular in Table~\ref{tab:sims-gene-expression} which shows that two
available programs, LISREL \citep{lisrel} and the R package `sem'
\citep{fox:2006}, fail to converge in a rather large number of problems.
This is in line with previous experience by other authors; see
e.g.~\cite{steiger:driving:2001}.  Our new RICF algorithm, on the other
hand, does not suffer from these problems.  It has clear convergence
properties (Theorem \ref{prop:converge}) and can handle problems with
several tens of variables (see Figure~\ref{fig:sim-model}).  In addition,
RICF has the desirable feature that it estimates parameters in closed form
(in a single cycle of its iterations) if this is possible. If applied to a
model based on a directed acyclic graph (DAG), the algorithm converges in a
single cycle and performs exactly the regressions commonly used for fitting
multivariate normal DAG models.  This feature and the fact that RICF can be
implemented using nothing but least squares computations make it an
attractive alternative to less specialized optimization methods.

\begin{figure}[tp]
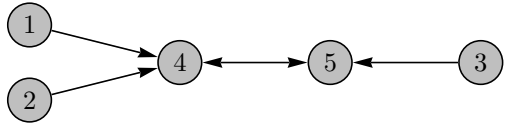

  \small
  \vspace{1.75cm}
  \begin{center}
  \begin{minipage}{7cm}
  \psset{linewidth=0.6pt,arrowscale=1.5 2,arrowinset=0.1}  
  \psset{fillcolor=lightgray, fillstyle=solid}          
  \newcommand{\myNode}[2]{\circlenode{#1}{\makebox[1.75ex]{#2}}}
  \rput(6.9, .5){\myNode{3}{$3$}}
  \rput(0.9, 0){\myNode{2}{$2$}} 
  \rput(0.9, 1){\myNode{1}{$1$}}
  \rput(4.9, .5){\myNode{5}{$5$}} 
  \rput(2.9, .5){\myNode{4}{$4$}}
  \ncline{->}{1}{4}
  \ncline{->}{2}{4}
  \ncline{->}{3}{5}
  \ncline{<->}{4}{5}
  \end{minipage}
  \caption{\label{fig:sur} Path diagram for seemingly unrelated
  regressions.} 
  \end{center}
\end{figure}

In another special case, namely seemingly unrelated regressions, RICF
reduces to the algorithm of \citet{telser:64}. A path diagram representing
seemingly unrelated regressions is shown in Figure \ref{fig:sur}.  The
variables $Y_1$, $Y_2$ and $Y_3$ are then commonly thought of as
covariates.  Since they have no spouses, the MLEs of the variances
$\omega_{11}$, $\omega_{22}$ and $\omega_{33}$ are equal to the empirical
variances $s_{11}$, $s_{22}$ and $s_{33}$.  For the remaining variables
$Y_i$, $i=4,5$, RICF performs regressions on both the ``covariates''
$Y_{\pa(i)}$ and the residual $\veps_j$, $j\in\{4,5\}\setminus\{i\}$.
These are precisely the steps performed by Telser.

Existing structural equation modelling software also fits models with
latent variables, whereas the RICF algorithm applies only to BAP models
without latent variables.  However, RICF could be used to implement the
M-step in the EM algorithm \citep{kiiveri:incomplete:1987} in order to fit
latent variable models.  This EM-RICF approach would yield an algorithm
with theoretical convergence properties.

Finally, we emphasize that the
RICF algorithm is determined by the path diagram.  However, different path
diagrams may induce the same statistical model; recall point
(\ref{item:equiv}) in \S\ref{sec:challenges} in the introduction.  This
model equivalence of path diagrams may be exploited to find a diagram for
which the running time of RICF is short.  Relevant graphical constructions
for this problem are described in \cite{drton:jmlr} and
\cite{ali:richardson:spirtes:zhang:towards:2005}.

\subsection*{Acknowledgements} 

This work was supported by the U.S.~National Science Foundation through
grant DMS-0505612, DMS-0505865, and the Institute for Mathematics and its
Applications, and by the U.S.~National Institutes for Health
(R01-HG2362-3).

\begin{appendix}

\section{Proofs}
\label{sec:app-proof-identify}

\begin{proof}[Proof of Proposition~\ref{prop:fisherinfo}]
  Let $\beta$ and $\omega$ be the vectors of unconstrained elements in $B$
  and $\Omega$, respectively.  The second derivatives of the log-likelihood
  function with respect to $\beta$ and $\omega$ are:
  \begin{align}
    \label{eq:2ndderivbb}   
    \frac{\partial^2\ell(B,\Omega)}{\partial\beta\,\partial\beta^t}
   &=-N\cdot P^t\big(S\otimes\Omega^{-1}\big)P,\\
    \label{eq:2ndderivbo}   
    \frac{\partial^2\ell(B,\Omega)}{\partial\beta\,\partial\omega^t}
    &=-N\cdot
    P^t\big[S(\idmat-B)^t\Omega^{-1}\otimes\Omega^{-1}\big]Q,\\
    \label{eq:2ndderivoo}
    \frac{\partial^2\ell(B,\Omega)}{\partial\omega\,\partial\omega^t}
    &=-\frac{N}{2}\,Q^t\big\{
    \big[\Omega^{-1}\otimes\Omega^{-1}(\idmat-B)S(\idmat-B)^t\Omega^{-1}\big]\\  
    \nonumber
    &\qquad\qquad+ \big[\Omega^{-1}(\idmat-B)S(\idmat-B)^t
    \Omega^{-1}\otimes\Omega^{-1}\big]\big\}Q. 
  \end{align}
  Replacing $S$ by $\Ew[S]=(\idmat-B)^{-1}\Omega(\idmat-B)^{-t}$ in
  (\ref{eq:2ndderivbb})-~(\ref{eq:2ndderivbo}) yields the claim.
\end{proof}


\begin{proof}[Proof of Proposition~\ref{prop:chaingraph}]
  If $G$ is a bi-directed chain graph, then the submatrix
  $B_{\tau_t,\tau_t}=0$ for all $t$, while for $s\not= t$ we have
  $\Omega_{\tau_s,\tau_t}=0$.  In this case the second derivative of the
  log-likelihood function with respect to $\beta_{ij}$ and $\omega_{kl}$ is
  equal to $\partial^2\ell(B,\Omega)/\partial
  \beta_{ij}\,\partial\omega_{kl}
  =[(\idmat-B)^{-1}]_{jl}\,(\Omega^{-1})_{ik}$.  Now
  $[(\idmat-B)^{-1}]_{jl}$ may only be non-zero if $j=l$ or $l$ is an
  ancestor of $j$, that is, if there exists a directed path $l\to
  j_1\to\dots\to j_m\to j$ in $G$.  On the other hand,
  $(\Omega^{-1})_{ik}=0$ whenever $i$ and $k$ are not in the same chain
  component.  Therefore, the second derivative in \eqref{eq:2ndderivbo} is
  equal to zero.
  
  Conversely, it follows that the second derivative in
  \eqref{eq:2ndderivbo} vanishes for all parameters only if the graph
  belongs to the class of bi-directed chain graphs.
\end{proof}

\end{appendix}

{\small
  \setlength\bibsep{0pt}
  \bibliographystyle{plainnat}
  \bibliography{sem}
}

\end{document}